\documentclass[onecolumn]{IEEEtran}

\ifCLASSINFOpdf
  \usepackage[pdftex]{graphicx}

\else

\fi

\usepackage{amsmath,amssymb, amsthm}
\usepackage{multirow}
\usepackage{authblk}
\usepackage{cite}
\usepackage{url}
\usepackage{pgf}
\usepackage{color}
\usepackage{todonotes}
\usepackage{booktabs} 
\usepackage{algorithm}
\usepackage{algorithmic}

\newcommand{\mb}[1]{\ensuremath{\boldsymbol{#1}}}
\newcommand{\comment}[1]{}
\begin{document}

\title{Adaptive Robust Optimization with Dynamic Uncertainty Sets for Multi-Period Economic Dispatch under Significant Wind}

\author{\'Alvaro~Lorca and Xu~Andy~Sun
\thanks{\'A. Lorca and X. A. Sun are with the Department of Industrial and Systems Engineering, Georgia Institute of Technology, Atlanta, GA 30332 USA}
}

\markboth{September 2014}%
{Shell \MakeLowercase{\textit{et al.}}: Bare Demo of IEEEtran.cls for Journals}

\maketitle

\IEEEpeerreviewmaketitle

\begin{abstract}
The exceptional benefits of wind power as an environmentally responsible renewable energy resource have led to an increasing penetration of wind energy in today's power systems. This trend has started to reshape the paradigms of power system operations, as dealing with uncertainty caused by the highly intermittent and uncertain wind power becomes a significant issue. Motivated by this, we present a new framework using adaptive robust optimization for the economic dispatch of power systems with high level of wind penetration. In particular, we propose an adaptive robust optimization model for multi-period economic dispatch, and introduce the concept of \textit{dynamic uncertainty sets} and methods to construct such sets to model temporal and spatial correlations of uncertainty. We also develop a simulation platform which combines the proposed robust economic dispatch model with statistical prediction tools in a rolling horizon framework. We have conducted extensive computational experiments on this platform using real wind data. The results are promising and demonstrate the benefits of our approach in terms of cost and reliability over existing robust optimization models as well as recent look-ahead dispatch models.
\end{abstract}


\allowdisplaybreaks

\section{Introduction}\label{sec:intro}
\IEEEPARstart{T}{he} exceptional benefits of wind power as an environmentally responsible energy resource have led to the rapid increase of wind energy in power systems all over the world. At the same time, wind energy possesses some characteristics  drastically different from conventional generating resources in terms of high stochasticity and intermittency in production output. Due to this, deep penetration of wind power will introduce significant uncertainty to the short-term and real-time operation of power systems, in particular, to the day-ahead unit commitment (UC) and the real-time economic dispatch (ED) procedures. If the uncertainty of such variable resources is not managed properly, the system operator may have to face severe operating conditions such as insufficient ramping capabilities from the conventional generating resources due to the sudden strong loss of wind power, complicated by other contingencies, load surge, and transmission congestions \cite{ERCOT2008wind}. These arising challenges call for new methods and models for power systems operation, and have attracted significant interests from both the electricity industry and academia.

The current UC and ED procedures rely on a combination of optimization tools and operational rules. The main optimization models used for UC and ED are deterministic models, where the uncertainties, such as demand, are assumed to take nominal forecast values. To deal with unexpected contingencies and sudden demand surge, the deterministic optimization model is complemented by operational rules that require extra generation resources, the so-called reserves, to stay available for quick response. The discrepancy between the forecast and realization of uncertainty has been relatively small in power systems composed of conventional load and supply. However, as observed in the recent experience, operating power systems with high penetration of variable resources, especially wind power, requires new methods to deal with uncertainty. See \cite{Xie2011wind} for an overview of the challenges of integrating wind in power systems from the perspective of UC, ED, frequency regulation and planning.

Facing these challenges, both industry and academia have devoted much effort to improving the current ED practice.
In particular, dynamic dispatch models with look-ahead capabilities have gained renewed interests. The basic ideas can be traced back to \cite{Bechert1972} and \cite{Ross1980}. Recent works have made significant advancement. \cite{Xie2011power} presents a look-ahead ED model with new statistical methods for wind forecast. The Midcontinent ISO has proposed look-ahead ED models with ramping products  \cite{navid2012market}. And \cite{ortega2009estimating} studies the selection of spinning reserve requirements under generation outages and forecast errors of demand and wind power. All these models can be characterized as deterministic ED models. Their simple optimization structure, improved performance, and closeness to the current operation make them appealing candidates to impact industry practice. This motivates the present paper to propose further advances and compare with these promising models.

Stochastic optimization has been a popular approach and extensively studied in the literature especially for the day-ahead unit commitment operation.
For example, \cite{Takriti1996StochUC} proposes one of the first stochastic UC models. \cite{Wu2007StochSCUC, Wang2008SCUC} propose security constrained UC models and consider stochastic wind power generation. \cite{bouffard2008stochastic} presents a short-term forward electricity market-clearing model under net load uncertainty, for the purpose of allowing high penetrations of wind power while keeping the system secure.  \cite{morales2009economic} deals with the selection of spinning and nonspinning reserves through a market-clearing model under stochastic wind power generation. \cite{Tuohy2009} presents a stochastic UC model for significant wind and shows the benefits of more frequent planning and over a deterministic approach. \cite{Papavasiliou2011reserve} studies reserve requirements for wind integration using a stochastic UC model.  \cite{papavasiliou2013multiarea} proposes multiarea stochastic UC models for high wind penetration. \cite{Wang2012chance} proposes a chance-constrained two-stage stochastic UC for systems with wind power uncertainty.

Regarding stochastic ED, the literature is much less extensive. \cite{keshmiri2010stoch} presents a stochastic programming model without recourse actions for a single-period ED problem.  \cite{lee2013afrequency} presents a stochastic model of a single-period ED problem under post-contingency frequency constraints.
\cite{zhang2013risk} presents a chance-constrained look-ahead ED model where the probability of incurring lost load is constrained and a sampling based scenario approximation approach is used for dealing with wind power randomness, however, transmission constraints are not considered in this work to ease computational burden. We would like to note that UC and ED have quite significant differences in decision structures and therefore modeling considerations: the UC has a relatively clear two-stage decision making structure, whereas for ED, the modeling choices are more diverse. Constructing a stochastic ED model with proper decision structure and desirable computational properties merits further research efforts. The literature in this respect still leaves much room for new contributions.



Recently, robust optimization has emerged as an alternative methodology for optimization under uncertainty \cite{RobustBook2009, Bertsimas2011review}. Robust optimization provides several features that are particularly appealing to applications in power systems. In particular, the robust optimization approach seeks to optimize system performance in a controlled manner against the worst-case scenario, which is indeed consistent with the philosophy of the current operational practice; robust optimization provides a data-driven way to model uncertainty, which scales well with the increasing dimension of data and is flexible and practical for many situations; robust optimization models are usually computationally tractable for large-scale systems.

Recent works have proposed robust optimization models for UC problems \cite{Street2011contingency, Guan2012pumphydro, Sun2013robustUC, Zeng2012robustUC, Zhao2013multistage}. \cite{Street2011contingency} provides a robust formulation for the contingency constrained UC problem. \cite{Guan2012pumphydro}, \cite{Sun2013robustUC}, \cite{Zeng2012robustUC} present two-stage adaptive robust models, with commitment decisions in the first stage and dispatch decisions in the second stage. In \cite{Sun2013robustUC} a two-stage robust UC model with security constraints is formulated and tested on the power system operated by ISO New England. \cite{Guan2012pumphydro} deals with a formulation including pumped storage hydro under wind power output uncertainty. Hybrid models and alternative objectives have also been explored to mitigate the conservativeness of the robust solution \cite{Zhao2013Unify, Jiang2013regret}. Efficient solution methods for the two-stage robust UC have been proposed \cite{Guan2012pumphydro, Sun2013robustUC, Zeng2013ccg, lee2013modeling}. Recently, \cite{lee2013modeling} presents acceleration techniques based on cutting planes and column generation for solving the two-stage robust UC problem under full transmission line constraints.


On the other hand, the benefits of robust optimization for the ED operation has not been fully explored. \cite{Zheng2012robustED} presents a two-stage robust ED model for a single-period regulation dispatch problem, where the first stage corresponds to dispatch and regulation capacity decisions, and the second stage corresponds to the dispatch of automatic generation control (AGC), after observing demand. \cite{Jabr2013affineOPF} recently proposes a robust optimal power flow model using affine policies for the AGC dispatch under renewable energy uncertainty. Affine policy is an approximation to the fully adaptive policy used in \cite{Zheng2012robustED}; however, as argued in \cite{Jabr2013affineOPF}, affine dependence on uncertainty may be a more suitable form for AGC dispatch. The work in \cite{xie2013shortterm} applies two advanced statistical methods for wind forecasting, and integrates these models with a robust look-ahead ED. However, their model is of a static robust nature, which lacks the adaptability of a two-stage robust model proposed here; their model also relies on the existing types of uncertainty sets, which will be significantly improved by a new type of uncertainty sets proposed in this paper.

If we try to summarize the above works, we can draw the following observations: 1) there is a great amount of interests to improve the ED practice; in particular, the recently developed look-ahead ED models have attracted considerable attention in both academia and industry; 2) the existing works on power system operation under uncertainty have focused on UC problems in a day-ahead operating environment, while both stochastic and robust ED models are relatively less explored; 3) the existing robust UC and ED models have used a similar type of uncertainty sets, which we call \emph{static} uncertainty sets, whereas it is important to start considering uncertainty sets that can capture the highly dynamical and correlated variable resources such as wind power.

In this paper, we propose new robust optimization models for system dispatch under high wind penetration. In particular, the contributions of our paper are summarized below:
\begin{enumerate}
\item We propose a two-stage adaptive robust optimization model for the multi-period ED, which has a different decision structure from the existing two-stage robust UC and robust ED models. The proposed robust ED model is designed for a rolling-horizon operational framework to model the real time ED process.
\item We introduce a new type of uncertainty sets, the \emph{dynamic uncertainty sets}, as a modeling technique to account for the dynamic relationship between uncertainties across decision stages. Such uncertainty sets explicitly model temporal and spatial correlations in variable sources. We also propose a data-driven approach to construct such dynamic uncertainty sets, which is simple to implement in practice.
\item We develop a comprehensive simulation platform, which integrates the proposed robust ED model with statistical procedures for constructing dynamic uncertainty sets using real-time data. Extensive experiments are performed on this platform.
\end{enumerate}

The paper is organized as follows. Section \ref{Section:UncSets} introduces dynamic uncertainty sets and discusses practical construction methods. Section \ref{Section:ModelAndMethod} proposes the adaptive robust multi-period ED model and solution methods. Section \ref{Section:EvaluationFramework} presents the simulation platform and the evaluation framework. Section \ref{Section:ComputationalExperiments} shows extensive computational experiments to demonstrate the effectiveness of our approach. Finally, Section \ref{Section:Conclusion} concludes.

\section{Dynamic uncertainty sets} \label{Section:UncSets}
In robust optimization, uncertainty is modeled through uncertainty sets, which are the building blocks of a robust optimization model and have direct impact on its performance. We may summarize three criteria for constructing uncertainty sets as follows. A well constructed uncertainty set should 1) capture the most significant aspects of the underlying uncertainty, 2) balance robustness and conservativeness of the robust solution, and 3) be computationally tractable.

\subsection{Static uncertainty sets}
Previous works on robust UC have focused on static uncertainty sets, and have treated uncertainty resources of different characteristics in an aggregated, indistinguishing way, see for example \cite{Sun2013robustUC, Guan2012pumphydro, Zeng2012robustUC}.
More specifically, consider the following uncertainty set for net demand vector $\mb{d}_{t} = (d_{1t},...,d_{N^dt})$:
\begin{align}
\notag & \mathcal{D}_t = \Bigg\{ \mb{d}_{t}: \; \sum_{j \in \mathcal{N}^d} \frac{|d_{jt}-\overline{d}_{jt}|}{\hat{d}_{jt}} \leq \Gamma^d \sqrt{N^d},\notag\\
& \qquad\qquad d_{jt} \in [\overline{d}_{jt} - \Gamma^d \hat{d}_{jt},\overline{d}_{jt} + \Gamma^d \hat{d}_{jt}] \; \forall \, j \in \mathcal{N}^d \Bigg\},\label{eq:sus}
\end{align}
where $\mathcal{N}^d, N^d$ denote the set and the number of loads, and $d_{jt}$ is the net demand of load $j$ at time $t$. According to \eqref{eq:sus}, $d_{jt}$ lies in an interval centered around the nominal value $\overline{d}_{jt}$ with a width determined by the deviation $\hat{d}_{jt}$. Further, the size of the uncertainty set is controled by $\Gamma^d$. If $\Gamma^d = 0$, $\mathcal{D}^t = \{\overline{\mb{d}}_{t}\}$, corresponding to a singleton set of the nominal demand. As $\Gamma^d$ increases, more demand vectors are contained in the uncertainty set, thus increasing the protection of the robust solution against larger demand variations.

The above uncertainty set is called \emph{static} uncertainty set, because the uncertainties at later time periods are \textit{independent} of those in earlier periods. That is, the dynamics of uncertainty evolution over time is not explicitly captured. Some recent work proposed additional budget constraints over time periods (e.g. \cite{Guan2012pumphydro,Zeng2012robustUC}). The modified uncertainty set imposes a coupling of uncertainty between time periods and uncertain sources, however, similar to \eqref{eq:sus}, it still does not directly characterize the temporal and spatial correlations of uncertainty; also, by coupling through the entire horizon, the realization of uncertainty breaks the time causality with past depending on the future realization. Yet another drawback of existing models is that uncertain sources of different nature are treated indistinguishably. For example, the uncertainty characteristics of wind power output are different from those of the conventional load, yet the existing proposals consider aggregated net load as the primitive uncertainty \cite{Sun2013robustUC, Guan2012pumphydro, Zeng2012robustUC}. Demand uncertainty is usually much less pronounced and less dynamic than wind, therefore, a static uncertainty set as \eqref{eq:sus} is an appropriate model. However, it is important to explore well suited uncertainty models for wind, specially for high level penetration of such variable resources.


\subsection{Dynamic uncertainty sets}
To explicitly model the correlation between multiple uncertain resources within one time period as well as the dynamics of each uncertain resource evolving over time periods, we propose the following general form of uncertainty sets, called \emph{dynamic uncertainty sets}: For each time $t$,
\begin{align}\label{eq:dus}
\Xi_t(\mb{\xi}_{[1:t-1]}) =\left\{\mb{\xi}_t : \exists \mb{u}_{[t]}\; \mbox{s.t.} \; f(\mb{\xi}_{[t]},\mb{u}_{[t]})\leq \mb{0} \right\},
\end{align}
where $\mb{\xi}_{[t_1:t_2]}\triangleq(\mb{\xi}_{t_1},\dots,\mb{\xi}_{t_2})$ and in shorthand $\mb{\xi}_{[t]}\triangleq \mb{\xi}_{[1:t]}$. In (\ref{eq:dus}), the uncertainty vector $\mb{\xi}_t$ explicitly depends on uncertainty at stages before time $t$ and the $\mb{u}$'s are auxiliary variables, $f(\mb{\xi}_{[t]},\mb{u}_{[t]})$ is a vector of convex functions that characterize the dynamics of uncertainty evolution. For the uncertainty set to be computationally tractable, $f$ should be semi-definite representable \cite{RobustBook2009}.

As an illustrative example, the dynamic uncertainty set could represent a dynamic interval for $\mb{\xi}_t$:
\begin{align*}
\mb{\xi}_t\in\left[\underline{\mb{\xi}}_t(\mb{\xi}_{[t-1]}),\overline{\mb{\xi}}_t(\mb{\xi}_{[t-1]})\right],
\end{align*}
where the upper and lower bounds of the interval at time $t$, namely $\underline{\mb{\xi}}_t(\mb{\xi}_{[t-1]})$ and $\overline{\mb{\xi}}_t(\mb{\xi}_{[t-1]})$, are functions of uncertainty realizations in previous time periods, rather than fixed values as in static uncertainty sets \eqref{eq:sus}.

A simple and useful specialization of (\ref{eq:dus}) is the linear dynamic uncertainty set, given as
\begin{align}\label{eq:lineardus}
\sum_{\tau=1}^t \left(\mb{A}_{\tau} \mb{\xi}_{\tau} + \mb{B}_{\tau}\mb{u}_{\tau}\right) \leq \mb{0},
\end{align}
which mimics linear dynamics and is also computationally tractable.
In the following, we will propose a specific method for constructing linear dynamic uncertainty sets using time series analysis tools.

\subsection{Constructing dynamic uncertainty sets for wind power}
The proposed dynamic uncertainty set (\ref{eq:dus}) is very general. In this section, we present a specific method to construct a dynamic uncertainty set for wind power using linear systems (\ref{eq:lineardus}). The key idea is to fuze time series models with the concept of dynamic uncertainty sets.

We denote the wind speed vector of multiple wind farms at time $t$ as $\mb{r}_t=(r_{1t},\dots,r_{N^wt})$, where $r_{it}$ is the wind speed at wind farm $i$ and time $t$.  Define the dynamic uncertainty set for $\mb{r}_t$ as:
\begin{subequations}\label{eq:windspeed_dus}
\begin{align}
\mathcal{R}_t({\mb{r}}_{[t-L:t-1]}) &=\Big\{\mb{r}_t : \;\exists\,\widetilde{\mb{r}}_{[t-L:t]},\; \mb{u}_t \; \quad \mbox{s.t.} \notag\\
\label{DefUnc}  & \quad \mb{r}_{\tau} = \mb{g}_{\tau} + \widetilde{\mb{r}}_{\tau}\;\;\,\forall \tau=t-L,\dots,t \\
\label{ARUnc}   & \quad \left.\widetilde{\mb{r}}_t = \sum_{s=1}^L \mb{A}_{s} \widetilde{\mb{r}}_{t-s} + \mb{B} \mb{u}_t  \right.\\
\label{uConst1} & \quad \sum_{i \in \mathcal{N}^w} |u_{it}| \leq \Gamma^w \sqrt{N^w}\\
\label{uConst2} & \quad |u_{it}| \leq \Gamma^w \quad\quad\; \forall i \in \mathcal{N}^w \\
\label{rConst}   & \quad  \mb{r}_t \geq \mb{0} \Big\},
\end{align}
\end{subequations}
where vectors $\mb{r}_{t-L},\dots,\mb{r}_{t-1}$ are the realizations of wind speeds in periods $t-L,\dots,t-1$. Eq. \eqref{DefUnc} decomposes wind speed vector $\mb{r}_{\tau}$ as the sum of a seasonal pattern $\mb{g}_{\tau}$, which is pre-estimated from wind data, and a residual component $\widetilde{\mb{r}}_{\tau}$ which is the deviation from $\mb{g}_{\tau}$. Eq. \eqref{ARUnc} is the key equation that represents a linear dynamic relationship involving the residual $\widetilde{\mb{r}}_t$ at time $t$, residuals realized in earlier periods $t-L$ to $t-1$, and an error term $\mb{u}_t$. The parameter $L$ sets the relevant time lags. In Eq. (\ref{ARUnc}), matrices $\mb{A}_s$'s capture the temporal correlation between $\mb{r}_t$ and $\mb{r}_{t-s}$, and $\mb{B}$ specifically captures the spatial relationship of wind speeds at adjacent wind farms at time $t$. Eq. \eqref{uConst1}-\eqref{uConst2} describe a budgeted uncertainty set for the error term $\mb{u}_t$, where $\Gamma^w$ controls its size, and (\ref{rConst}) avoids negative wind speeds. $\mathcal{N}^w$ and $N^w$ denote the set and number of wind farms, respectively.



Using the above uncertainty sets \eqref{eq:windspeed_dus} for wind speeds, we can further construct dynamic uncertainty sets for wind power through power curve approximations. In particular, we denote the \emph{available wind power} of wind farm $i$ at time $t$ as $\overline{p}_{it}^w$. Given the wind speed $r_{it}$, $\overline{p}_{it}^w$ is described by the following constraints
\begin{align}
\label{eq:PowerCurve} \overline{p}_{it}^w \geq h_{ik}^0 + h_{ik} \, r_{it} \quad \forall i \in \mathcal{N}^g, \, k=1,\dots,K,
\end{align}
where parameters $h_{ik}^0, h_{ik}$ are determined based on a convex piecewise linear approximation with $K$ pieces of the increasing part of the power curve at wind farm $i$ (in our experiments, we use the power curve of GE 1.5MW wind turbine to approximate the aggregated output of a wind farm). Although (\ref{eq:PowerCurve}) allows available wind power to exceed $\max_k \{h_{ik}^0 + h_{ik} r_{it}\}$, the robust optimization model described in Section \ref{Section:ModelAndMethod} will always ensure that the available wind power lies on the power curve including the plateau part for wind speed exceeding a cut-off value.

The dynamic uncertainty set of the available wind power $\overline{\mb{p}}_{t}^w$ is thus defined as
\begin{align}\label{eq:windpower_dus}
\overline{\mathcal{P}}^{w}_t(\mb{r}_{[t-L:t-1]}) = & \Big\{\overline{\mb{p}}_{t}^w :\; \exists \mb{r}_t \in \mathcal{R}_t({\mb{r}}_{[t-L:t-1]}) \notag \\
&\hspace{2cm} \text{s.t. \eqref{eq:PowerCurve} is satisfied}\Big\},
\end{align}
based on which we can define the uncertainty set for the trajectory of available wind power in time periods $2$ through $T$, namely $\overline{\mb{p}}^w = (\overline{\mb{p}}^w_2, \dots, \overline{\mb{p}}^w_T)$, as
\begin{align}\label{eq:windpowertrajectory_dus}
\overline{\mathcal{P}}^{w} = & \Big\{(\overline{\mb{p}}^w_2, \dots, \overline{\mb{p}}^w_T) : \; \exists (\mb{r}_2, \dots, \mb{r}_T) \; \mbox{s.t.} \; \mb{r}_t \in \mathcal{R}_t({\mb{r}}_{[t-L:t-1]})\notag\\
&\quad\mbox{and} \; \text{(\ref{eq:PowerCurve}) is satisfied for $t=2,\dots,T$} \Big\},
\end{align}
which is used in the robust ED model.

As a summary, we propose dynamic uncertainty sets \eqref{eq:windspeed_dus} and \eqref{eq:windpower_dus} to capture the intrinsic temporal dynamics and spatial correlations of the wind power. We also distinguish wind power uncertainty from conventional demand uncertainty, which is modeled by traditional uncertainty sets \eqref{eq:sus}. The proposed dynamic uncertainty set formulation (\ref{eq:dus}) is quite general. The specific models for wind speed (\ref{eq:windspeed_dus}) and wind power (\ref{eq:windpower_dus}) present one example for its implementation. Other models may be constructed using more sophisticated statistical tools. For example,  the coefficient matrices $\mb{A}_s$ and $\mb{B}$ can be made time dependent as $\mb{A}_{st}$ and $\mb{B}_t$ using dynamic vector autoregression methods. Also, $\mb{r}_t$ can be replaced by a nonlinear transformation of wind speed to improve estimation accuracy.
However, there is always a tradeoff between model complexity and performance. Our experiments show the above simple models \eqref{eq:windspeed_dus}-\eqref{eq:windpowertrajectory_dus} achieve a substantial improvement over existing static uncertainty sets. See more discussion in Section \ref{Sec:ParameterEstimation} for parameter estimation and possible extensions for the dynamic uncertainty sets.

\section{Adaptive robust economic dispatch formulation and solution method} \label{Section:ModelAndMethod}
\subsection{Mathematical formulation}
In this section, we propose an adaptive robust optimization model for the multi-period ED problem. In particular, the ED problem with $T$ time periods is formulated as a two-stage adaptive robust model in the following way. The first-stage of the robust ED model comprises the current time period $t=1$, while the second-stage comprises future time periods $t=2,...,T$. In the first-stage, the decision maker observes demand and available wind power at the current time period, and determines the dispatch solution, which will be implemented right away for time $t=1$. Given the first-stage decision, the second-stage of the robust ED model computes the worst-case dispatch cost for the remaining time periods in the dispatch horizon. The overall robust ED model minimizes the total cost of dispatch at the current time period and the worst-case dispatch cost over the future periods.

We denote $\mb{x} = (\mb{p}^g_{1}, \mb{p}^w_{1})$ as the vector of first-stage dispatch decisions, composed of generation of thermal units ($\mb{p}^g_{1}$) and wind farms ($\mb{p}^w_{1}$). Note that we allow wind generation to be dispatchable. The uncertainty includes both conventional load $\mb{d}=(\mb{d}_2,\dots,\mb{d}_T) \in \mathcal{D}$ described by (\ref{eq:sus}) and the available wind power $\overline{\mb{p}}^w=(\overline{\mb{p}}^w_2, \dots, \overline{\mb{p}}^w_T) \in \overline{\mathcal{P}}^w$ described by the dynamic uncertainty set \eqref{eq:windpowertrajectory_dus}. The second-stage decisions are dispatch solutions  $\mb{y}=\left(\mb{p}^g_{t},\mb{p}^w_{t}, \forall t=2,\dots, T\right)$.

Mathematically, the two-stage robust multi-period ED model is formulated as follows,
\begin{align}
\label{Formulation} & \min\limits_{\mb{x} \in {\Omega}_1^{det}} \Bigg\{ \mb{c}^{\top} \mb{x} + \max\limits_{\mb{d} \in \mathcal{D}, \overline{\mb{p}}^w \in \overline{\mathcal{P}}^w} \; \min\limits_{\mb{y} \in {\Omega}(\mb{x},\mb{d},\overline{\mb{p}}^w)} \; \mb{b}^{\top}\mb{y} \Bigg\},
\end{align}
where the first and second-stage costs are defined as
\begin{align*}
& \mb{c}^{\top} \mb{x} = \sum_{i \in \mathcal{N}^g} C^g_i p^g_{i1} + \sum_{i \in \mathcal{N}^w} C^w_i p^w_{i1} \\
& \mb{b}^{\top}\mb{y} = \sum_{t=2}^{T} \left(\sum_{i \in \mathcal{N}^g} C^g_i p^g_{it} + \sum_{i \in \mathcal{N}^w} C^w_i p^w_{it} \right),
\end{align*}
where $\mathcal{N}^g$ denotes the set of generators, and $C^g_i, C^w_i$ denote the variable costs of thermal generators and wind farms. We use linear dispatch costs, but it is straightforward to extend to piecewise linear approximations of nonlinear cost functions.

The feasible region $\Omega_1^{det}$ of the first-stage decision variables corresponds to the constraints of a single-period dispatch problem, that is
\begin{subequations}\label{eq:1ststage_constr}
\begin{align}
\Omega_1^{det} = \Bigg\{&\mb{x} = (\mb{p}_1^g,\mb{p}_1^w):\;
\label{LimitsG1} \underline{p}_{i1}^g \leq p^g_{i1} \leq \overline{p}_{i1}^g \quad\forall \, i \in \mathcal{N}^g \\
\label{LimitsW1} & 0 \leq p^w_{i1} \leq p_i^{w,max}\quad \forall \, i \in \mathcal{N}^w \\
\label{AvailableWind1} & p^w_{i1} \leq \overline{p}_{i1}^{w,det} \quad \forall \, i \in \mathcal{N}^w \\
\label{RampG1} & -RD^g_i \leq p^g_{i1} - p^g_{i0} \leq RU^g_i \quad \forall \, i \in \mathcal{N}^g \\
\label{RampW1} & -RD^w_i \leq p^w_{i1} - p^w_{i0} \leq RU^w_i \quad \forall \, i \in \mathcal{N}^w \\
\label{LineFlow1} & \hspace{-10mm}\left|\mb{\alpha}_l^{\top} (\mb{E}^g \mb{p}^g_{1} + \mb{E}^w \mb{p}^w_{1} - \mb{E}^d \mb{d}^{det}_{1})\right| \leq f_l^{max} \;\forall \, l \in \mathcal{N}^l \\
\label{Balance1} & \sum_{i \in \mathcal{N}^g} p^g_{i1} + \sum_{i \in \mathcal{N}^w} p^w_{i1} = \sum_{j \in \mathcal{N}^d} d^{det}_{j1}\Bigg\},
\end{align}
\end{subequations}
where $\underline{p}_{it}^g, \overline{p}_{it}^g$ are the minimum and maximum power outputs of thermal generator $i$ at time $t$; $p_i^{w,max}$ is the maximum power output at wind farm $i$, representing the cut-off level of the power curve;  $\overline{p}_{i1}^{w,det}$ denotes the available wind power of wind farm $i$ observed at current time $t=1$; $RD^g_i, RU^g_i$ are the ramp-down and ramp-up rates of thermal generators (similarly, $RD^w_i, RU^w_i$ for wind farms); $\mathcal{N}^l$ is the set of transmission lines; $\mb{\alpha}_l$ is the network shift factor for line $l$; $\mb{E}^d, \mb{E}^g, \mb{E}^w$ are the network incidence matrices for loads, thermal generators and wind farms; $f_l^{max}$ is the flow limit on line $l$; $d^{det}_{j1}$ denotes the observed electricity demand at load $j$ and time $t=1$. Constraints \eqref{LimitsG1}, \eqref{LimitsW1} and \eqref{AvailableWind1} enforce generation limits for thermal generators and wind farms, with \eqref{AvailableWind1} ensuring that generation of wind farms does not exceed the available wind power at time $t=1$. \eqref{RampG1} and \eqref{RampW1} enforce ramping rate limits for thermal generators and wind farms. \eqref{LineFlow1} represents line flow limits. \eqref{Balance1} represents energy balance.

Constraints in the second-stage problem  are parameterized by the first-stage decision variables and uncertain parameters realized in the uncertainty sets. The feasible region of the second-stage dispatch decison $\mb{y}=\left(\mb{p}^g_{t},\mb{p}^w_{t}, \forall t=2,\dots, T\right)$ is defined as
\begin{subequations}\label{eq:2ndstage_constr}
\begin{align}
\Omega(\mb{x},\mb{d},&\overline{\mb{p}}^w)=\Bigg\{\mb{y}:\; \; \mbox{s.t.} \;\;\forall t=2,\dots,T \notag\\
\label{LimitsG2toT} &  \underline{p}_{it}^g \leq p^g_{it} \leq \overline{p}_{it}^g \quad \forall \, i \in \mathcal{N}^g,  \\
\label{Limitswmax2toT} & 0 \leq p^w_{it} \leq p_i^{w,max} \quad \forall \, i \in \mathcal{N}^w, \\
\label{AvailableWind2toT} & p^w_{it} \leq \overline{p}_{it}^w \quad \forall \, i \in \mathcal{N}^w,   \\
\label{RampG2toT} & -RD^g_i \leq p^g_{it} - p^g_{i,t-1} \leq RU^g_i \quad \forall \, i \in \mathcal{N}^g,  \\
\label{RampW2toT} & -RD^w_i \leq p^w_{it} - p^w_{i,t-1} \leq RU^w_i \quad \forall \, i \in \mathcal{N}^w, \\
\label{LineFlow2toT} &\hspace{-10mm}\left|\mb{\alpha}_l^{\top} (\mb{E}^g \mb{p}^g_{t} + \mb{E}^w \mb{p}^w_{t} - \mb{E}^d \mb{d}_{t})\right| \leq f_l^{max} \;\forall \, l \in \mathcal{N}^l  \\
\label{Balance2toT} & \sum_{i \in \mathcal{N}^g} p^g_{it} + \sum_{i \in \mathcal{N}^w} p^w_{it}= \sum_{j \in \mathcal{N}^d} d_{jt}\Bigg\},
\end{align}
\end{subequations}
where \eqref{LimitsG2toT}-\eqref{Balance2toT} are similar constraints as in (\ref{eq:1ststage_constr}), except that they are enforced for each time period $t=2,\dots, T$. Notice that \eqref{Limitswmax2toT}-\eqref{AvailableWind2toT} ensure that the dispatched wind generation is upper bounded by the minimum between the cut-off level $p_i^{w,max}$ and the available wind power $\overline{p}_{it}^w$. Also note that the first-stage dispatch decision is involved in constraints \eqref{RampG2toT}-\eqref{RampW2toT} to satisfy ramping constraints.

A few remarks are in order. First, \eqref{Formulation} is a fully adaptive robust optimization model, namely the second-stage dispatch decision adapts to every realization of the uncertainty in the best possible way, which is similar to the existing robust UC model proposed in \cite{Sun2013robustUC}. Second, there is  a key difference between the two-stage structure of the proposed robust ED (\ref{Formulation}) and the existing two-stage robust UC models. In particular, the decision stages of (\ref{Formulation}) correspond to the actual time periods, so that the first-stage decision can be directly used in the dispatch at the current period, and the dispatch decisions in the second stage can be re-optimized in the following periods. In comparison, the two-stage robust UC models have UC decisions in the first stage and dispatch decisions in the second stage, both for the \emph{entire horizon}. Third, the two-stage structure of the robust ED model makes it convenient to incorporate into the real-time dispatch procedure. In particular, the robust ED model can be implemented in a rolling horizon framework; the dynamic uncertainty sets can also be updated periodically when new information is available. Fourth, the use of the DC power flow is consistent with the industry practice \cite{FERC2011ACOPF} and recent works in robust ED \cite{Zheng2012robustED, Jabr2013affineOPF}. AC power flow feasibility can be enforced by introducing an AC power flow module.  Thus, to emphasize the key proposal of the paper, we keep with the simple DC power flow model. Fifth, the robust ED model can also readily include convex piecewise linear costs.

\subsection{Solution method} \label{Subsection:SolutionMethod}
Several methods have been reported in the literature for solving two-stage adaptive robust optimization problems \cite{Sun2013robustUC,Guan2012pumphydro,Zeng2013ccg}. In \cite{Guan2012pumphydro}, a Benders decomposition approach is proposed to solve the outer level problem and an exact method for the second-stage problem. In \cite{Zeng2013ccg}, a constraint and column generation (C\&CG) technique is proposed and rigorously analyzed; an exact method using mixed-integer reformulations is proposed for the second-stage problem. In \cite{Sun2013robustUC}, a modified Benders decomposition framework is proposed for the outer level problem with an efficient heuristic method for the second-stage problem. The key modification to the traditional Benders decomposition is to add the worst-case extreme point and the associated dispatch constraints to the outer level problem in each iteration (see \cite[Section IV]{Sun2013robustUC}). This is similar to the idea behind constraint and column generation in \cite{Zeng2013ccg}.

Problem \eqref{Formulation} can be equivalently stated as:
\begin{align} \label{eq:Formulation2}
\min\limits_{\mb{x},\eta} \left\{ \mb{c}^\top \mb{x} + \eta: \; \eta \geq Q(\mb{x}), \; \mb{x} \in \Omega_1^{det} \right\},
\end{align}
with
\begin{align} \label{eq:Qbeforedual}
Q(\mb{x}) = \max\limits_{\mb{\xi} \in \Xi} \; \min\limits_{\{\mb{y}: \;\; \mb{G} \mb{y} \geq \mb{h} - \mb{E} \mb{x} - \mb{M} \mb{\xi}\}} \;\; \mb{b}^\top \mb{y},
\end{align}
where $\mb{\xi} = (\mb{d},\overline{\mb{p}}^w)$, $\Xi = \mathcal{D} \times \overline{\mathcal{P}}^w$, and the feasible region $\{ \mb{y}: \mb{G} \mb{y} \geq \mb{h} - \mb{E} \mb{x} - \mb{M} \mb{\xi} \}$ represents the dispatch constraints in \eqref{eq:2ndstage_constr}. Problem \eqref{eq:Formulation2} is equivalent to:
\begin{subequations}\label{eq:ccgformulation}
\begin{align}
\label{eq:ccgformulationobjective} \min\limits_{\mb{x} \in \Omega_1^{det},\, \eta, \, \{\mb{y}_l\}} & \quad \mb{c}^\top \mb{x} + \eta \\
\label{eq:ccgformulationopt} \mbox{s.t.} & \quad \eta \geq \mb{b}^\top \mb{y}_l \quad \forall \; l \\
\label{eq:ccgformulationfeas} & \quad \mb{E}\mb{x} + \mb{G} \mb{y}_l \geq \mb{h} - \mb{M} \mb{\xi}_l^* \quad \forall \; l,
\end{align}
\end{subequations}
where $\{\mb{\xi}_l^*\}_{l=1}^M$ is the set of extreme points of $\Xi$, and for each $l$, $\mb{y}_l$ is a vector of second-stage decisions associated to $\mb{\xi}_l^*$.  \eqref{eq:ccgformulation} is the outer level problem, which shows a nice structure suitable for constraint generation. Indeed, \eqref{eq:ccgformulation} can be efficiently solved by adding $(\mb{\xi}_l^*, \mb{y}_l)$ and the associated constraints iteratively \cite{Sun2013robustUC, Zeng2013ccg}.

In every iteration of this algorithm, $Q(\mb{x})$ must be evaluated, which involves solving a nonconvex max-min problem. Previous work has dealt with this problem using outer-approximation techniques \cite{Sun2013robustUC} and exact methods based on mixed-integer programming (MIP) reformulations \cite{Guan2012pumphydro,Zeng2012robustUC,Zeng2013ccg}. As will be demonstrated in the computational experiments (Section \ref{subsection:alternatingmethodperformance}), the MIP method is time consuming for solving \eqref{eq:Qbeforedual}. Instead, we apply a simple ``alternating direction algorithm'' \cite{konno1976cutting}. Taking the dual over the inner $\min$ in \eqref{eq:Qbeforedual} we obtain
\begin{equation} \label{eq:Qafterdual}
Q(\mb{x}) = \;  \max\limits_{\mb{\xi} \in \Xi, \mb{\pi} \in \Pi} \quad \mb{\pi}^\top (\mb{h}-\mb{E}\mb{x}-\mb{M}\mb{\xi}),
\end{equation}
where $\Pi = \{\mb{\pi} \geq \mb{0}: \mb{\pi}^\top \mb{G} = \mb{b}\}$. For this \emph{bilinear} program with separate polyhedral feasible regions $\Xi$ and $\Pi$, the alternating direction algorithm optimizes over $\mb{\pi}$ with $\mb{\xi}$ fixed, then over $\mb{\xi}$ with $\mb{\pi}$ fixed, and alternates; each of these iterations solves a linear program which achieves the optimum at an extreme point of the corresponding polyhedron $\Xi$ or $\Pi$. The alternating algorithm is formally presented below.


\begin{algorithm}\label{algorithm:AlternatingAlgorithm}
\caption{Alternating Direction (AD) algorithm}
\small
\begin{algorithmic}[1]
\STATE{Start with some $\mb{\xi}' \in \Xi$}
\REPEAT
    \STATE{Solve $(*)$: $C \leftarrow \max_{\mb{\pi} \in \Pi} \; \mb{\pi}^\top (\mb{h}-\mb{E}\mb{x}-\mb{M}\mb{\xi}')$}
    \IF{$C < \infty$}
        \STATE{Let $\mb{\pi}'$ be an optimal solution of $(*)$}
        \STATE{Solve $C' \leftarrow \max_{\mb{\xi} \in \Xi} \; {\mb{\pi}'}^\top (\mb{h}-\mb{E}\mb{x}-\mb{M}\mb{\xi})$ and let $\mb{\xi}'$ be its optimal solution}
    \ELSE
        \STATE{$C' \leftarrow \infty$}
    \ENDIF
\UNTIL{$C' = \infty$ \OR $C' - C \leq \delta$}
\STATE{\textbf{output}: $C'$ as estimate of $Q(\mb{x})$ with solution $\mb{\xi}'$}
\end{algorithmic}
\end{algorithm}

This alternating direction method always converges to a KKT point of \eqref{eq:Qafterdual}. The proof is omitted to save space. Section \ref{subsection:alternatingmethodperformance} also shows empirical evidence that this heuristic achieves good solution quality and fast convergence on the second-stage problem, comparing to the MIP method.

\begin{figure}[h!]
\centering
\includegraphics[width=3.5in]{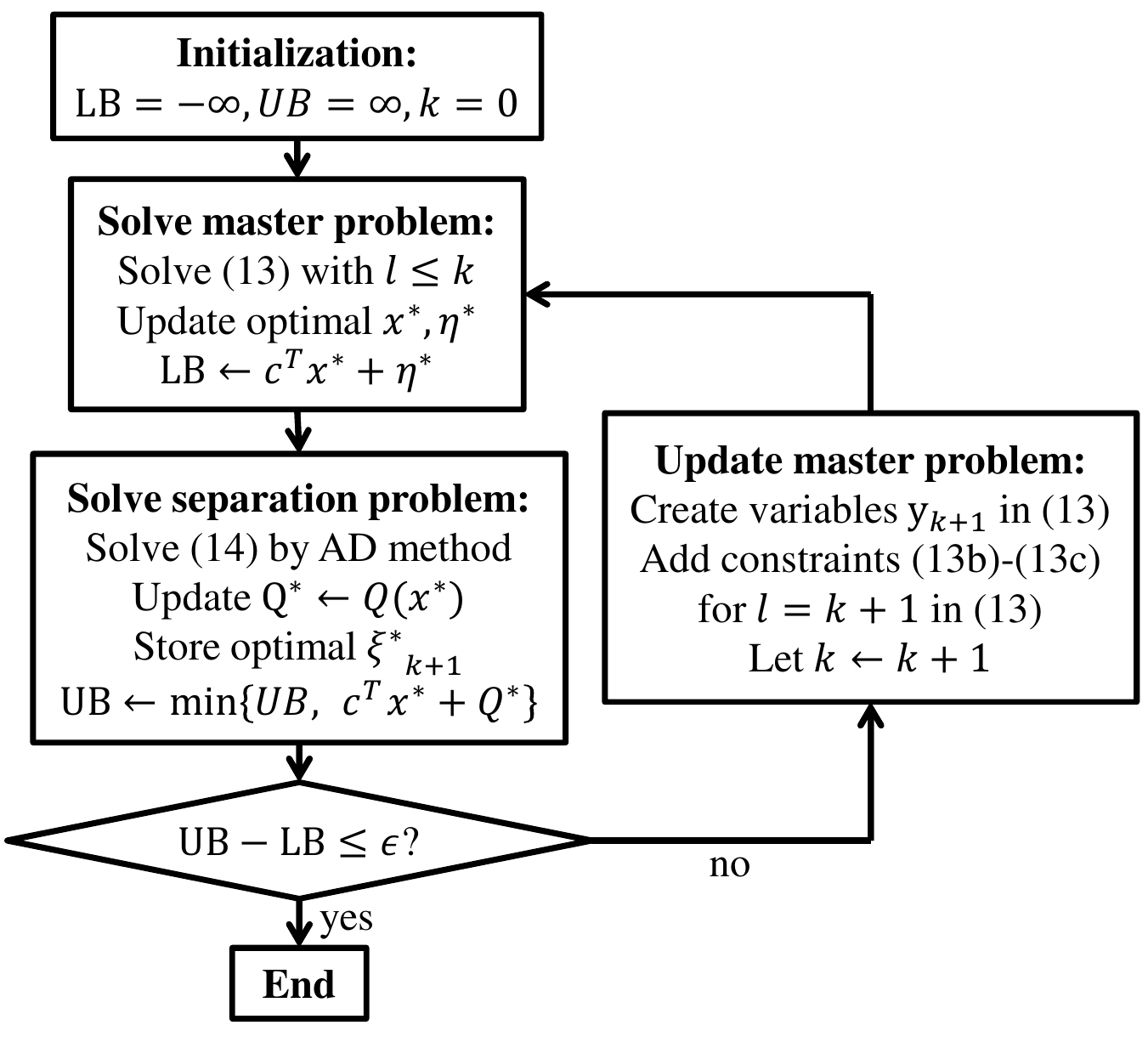}
\caption{Flow chart for the overall two-level algorithm.}
\label{Figure:AlgorithmFlowChart}
\end{figure}

The overall two-level algorithm is presented in Fig. \ref{Figure:AlgorithmFlowChart}.

\section{Simulation platform and evaluation metrics} \label{Section:EvaluationFramework}
\begin{figure}[h!]
\centering
\includegraphics[width=3.5in]{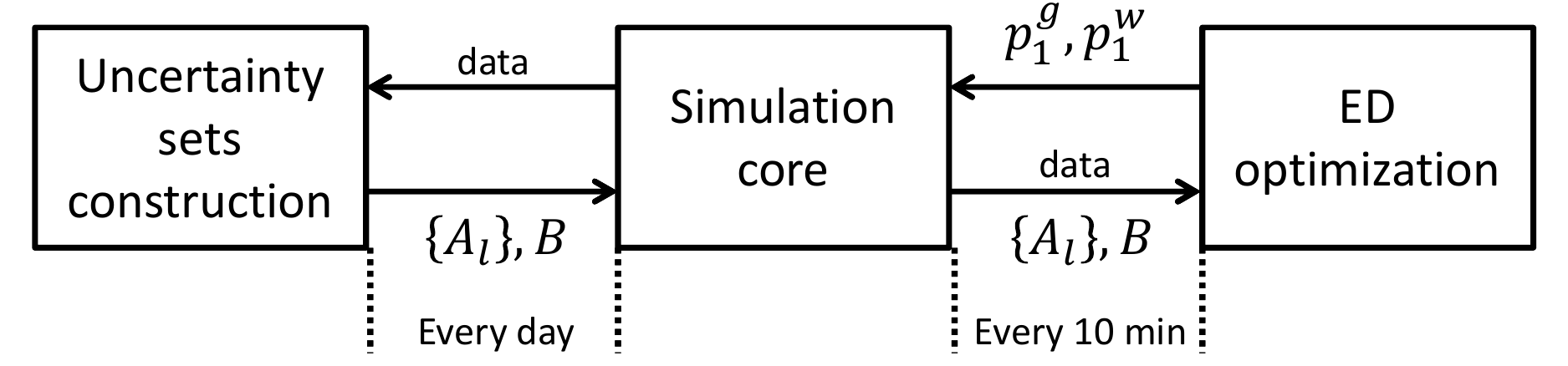}
\caption{Simulation platform integrating ED optimization engine and data analysis tools for uncertainty model construction.}
\label{Figure:platform}
\end{figure}

In this Section, we describe the simulation platform and evaluation metrics for the proposed robust model. The motivation is to have a realistic
simulation environment that integrates the dispatch optimization model with data analysis procedures which dynamically update the parameters in the optimization and uncertainty models. Fig. \ref{Figure:platform} illustrates the simulation process.

The simulation process is implemented in a rolling horizon framework. At each time period, the robust ED model is solved over a time window of $T$ time periods. The first-stage dispatch solution for the current time period is implemented, while the second-stage dispatch solutions for remaining periods are not materialized; the time horizon rolls forward by one time interval, where new realizations of demand and available wind power are observed, and dynamic uncertainty sets are periodically re-estimated and updated with the new observations (see Section \ref{Sec:ParameterEstimation}). In order to focus the comparison on the ED policies, the simulation process uses a simplified UC schedule where all thermal generators are on all the time. In the future, we would like to extend the simulation framework to integrate UC decisions into the policy evaluation.

\begin{figure}[h!]
\centering
\includegraphics[width=3in]{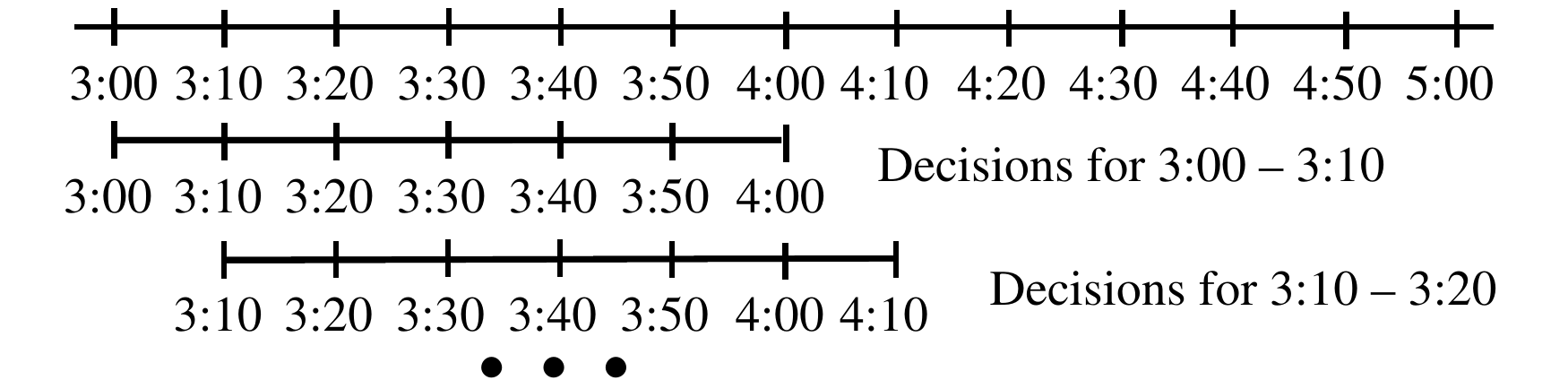}
\caption{Concept of rolling horizon with 10 minute time periods and $T=6$}
\label{Figure:rollinghorizon}
\end{figure}

We compare different ED models by evaluating the average and standard deviation (std) of the production cost for every 10 minutes dispatch interval, which includes both generation cost and penalty cost resulting from the use of expensive fast-start units or load shedding.

\subsection{Estimating the parameters of the dynamic uncertainty set for wind speeds} \label{Sec:ParameterEstimation}
In order to estimate the parameters of model \eqref{eq:windspeed_dus}, consider the following time series model:
\begin{subequations}\label{eq:windspeed_timeseriesmodel}
\begin{align}
\label{DefTimeSeriesModel} & \quad \mb{r}_{t} = \mb{g}_{t} + \widetilde{\mb{r}}_{t} \quad \forall t \\
\label{ARTimeSeriesModel} & \quad \left.\widetilde{\mb{r}}_t = \sum_{s=1}^L \mb{A}_{s} \widetilde{\mb{r}}_{t-s} + \mb{\epsilon}_t \right. \quad \forall t,
\end{align}
\end{subequations}
where $\mb{r}_{t}$ is the vector of wind speeds at time $t$, $\mb{g}_{t}$ corresponds to a deterministic seasonal pattern, and $\widetilde{\mb{r}}_{t}$ corresponds to the deviation of $\mb{r}_{t}$ from $\mb{g}_{t}$. In this model, $\widetilde{\mb{r}}_t$ follows a multivariate autoregressive process of order $L$, determined by the innovation process $\{\mb{\epsilon}_t\}$, where $\mb{\epsilon}_t$ is a vector of  normal random variables with mean \mb{0} and covariance matrix $\mb{\Sigma}$, and vectors $\mb{\epsilon}_t$ are independent across different time periods.

Once seasonal patterns have been identified, parameter $\mb{g}_t$ can be determined. For example, daily and semi-daily seasonalities could be used. In such a case, using a 10 min time interval we would have $g_{it} = a_i + b_i \cos(\frac{2 \pi t}{24 \times 6}) + c_i \sin(\frac{2 \pi t}{24 \times 6}) + d_i \cos(\frac{2 \pi t}{12 \times 6}) + e_i \sin(\frac{2 \pi t}{12 \times 6})$ (since $24 \times 6$ is the number of time periods in a day). Parameters $a_i,b_i,c_i,d_i,e_i$ can be estimated using linear regression \cite{Xie2011power}.

The parameters of the autoregressive component $\widetilde{\mb{r}}$, namely the matrices $\mb{A}_{s}$ and $\mb{\Sigma}$, can be estimated using statistical inference techniques developed for time series \cite{reinsel1997elementsbook}, for which many computational packages are available. $\mb{B}$ in \eqref{eq:windspeed_dus} is obtained from the Cholesky decomposition of $\mb{\Sigma}$.

The linear dynamic model \eqref{eq:windspeed_dus} and the associated estimation method are appealing in their simple structure, which serves well our goal to demonstrate the concept of dynamic uncertainty sets. Computational results also confirm their promising performance. Meanwhile, it is worth noting that the framework of dynamic uncertainty sets is flexible enough to incorporate more sophisticated statistical models, such as the ones proposed in  \cite{morales2010methodology}, where autoregressive processes are fitted to nonlinearly transformed wind speeds. Using a piecewise linear approximation similar to the one proposed in \eqref{eq:PowerCurve}, but this time for the transformed wind speed and wind power output, a dynamic uncertainty set can be again constructed using linear constraints.

\section{Computational experiments} \label{Section:ComputationalExperiments}
We conduct extensive computational experiments on the simulation platform to compare the proposed robust ED model and dynamic uncertainty sets with existing robust and deterministic dispatch models. The experiments are performed on the 14-bus and 118-bus IEEE test systems, both of which are modified to incorporate significant wind penetration. In the following, we introduce the detailed data for the 14-bus system, and present test results in Sections \ref{subsection:robustvsdeterministic} to \ref{subsection:alternatingmethodperformance}. The test results on the 118-bus system give a similar picture as the 14-bus system. The details are given in Section \ref{subsection:118bus}.

Table \ref{14busGenerators} summarizes $P_{min}, P_{max}$, 10-min ramping rates, and production costs of all three generators in the 14-bus system. The total generation capacity is 500MW. The system has 20 transmission lines and 11 conventional loads. The daily system demand is between 132.6MW and 319.1MW with an average of 252.5MW. The system has 4 wind farms, each with a capacity of 75MW (equivalent of 50 units of GE 1.5MW wind turbines). The total power output at each wind farm is approximated by a piecewise linear function of wind speed using the power curve data \cite{GEonepointfiveunit}. 

\begin{table}[h!]
\renewcommand{\arraystretch}{1}
\caption{Thermal generators in 14-bus system}
\label{14busGenerators}
\centering
\begin{tabular}{c|c|c|c|c}
\toprule
\textbf{Gen} & \textbf{Pmax} & \textbf{Pmin} & \textbf{Ramp} & \textbf{Cost} \\
 & \textbf{(MW)} & \textbf{(MW)} & \textbf{(MW/10min)} & \textbf{(\$/MWh)} \\
\midrule
1 & 300 & 50 & 5 & 20 \\
2 & 100 & 10 & 10 & 40 \\
3 & 100 & 10 & 15 & 60 \\
\bottomrule
\end{tabular}
\end{table}

The wind speed data is obtained from \cite{ceerewebsite} for four geographically adjacent locations with a 10-minute data interval. The average wind speeds at the four wind farms are 4.8, 5.6, 5.1, 5.5 m/s, respectively. Using the power curve, the average total available wind power is 104.2 MW, equivalent to a $34.7\%$ capacity factor, which is about $32.7\%$ of peak demand and $20\%$ of conventional generation capacity, representing a realistically high level of wind penetration. After removing stationary components, wind speeds at different sites present strong auto and cross correlation at several lags, which implies that the temporal and spatial dependencies are significant.  

The proposed robust ED model has $9$ time periods with a 10-min interval for each period (i.e. 1.5-hour look ahead). The robust ED model is evaluated on the simulation platform in the rolling-horizon framework. In particular, it is solved every 10 minutes over 35 days, for which real wind data is used for all wind farms. On each of the 35 days, the simulation engine updates the parameters of the dynamic uncertainty sets  (\ref{eq:windspeed_dus}) using the available wind data up to that day. The penalty cost is $C^+ = 6000$ \$/MWh for under-generation, and $C^- = 600$ \$/MWh for over-generation \cite{papavasiliou2013multiarea, Zhao2013Unify}.


The simulation platform is implemented in a Python environment, interfaced with Cplex 12.5.
Each robust ED takes less than a second to solve, and the entire simulation of 5040 periods takes about 40 minutes on a PC laptop with an Intel Core i3 at 2.1 GHz and 4GB memory.

Before presenting details, we first give a summary of the experiments and main results.  We compare the proposed robust ED with dynamic uncertainty sets versus (1) deterministic look-ahead dispatch and its variant with reserve rules; (2) robust dispatch with static uncertainty sets. The experiments show that adaptive robust ED with dynamic uncertainty sets significantly outperforms both alternative models by substantially reducing average production cost, the variability of the costs, and the probability of shortage events. Our experiments also show that the robust ED provides a \emph{Pareto frontier} for the tradeoff between cost and reliability, which provides an informative guideline for choosing uncertainty set parameters and system operating points.

\subsection{Robust ED versus look-ahead ED} \label{subsection:robustvsdeterministic}
In this section, we compare the proposed adaptive robust ED (Rob-ED) with the deterministic look-ahead dispatch (LA-ED). The robust ED model uses dynamic uncertainty sets \eqref{eq:windspeed_dus} and \eqref{eq:windpower_dus} with 6 time lags i.e. $L=6$. \comment{We choose $L=6$ because it outperforms other lag numbers.} The parameter $\Gamma^w$ controls the size of the uncertainty sets. Notice that when $\Gamma^w = 0$, the uncertainty set contains only one path of the forecasted wind speeds, the robust ED thus reduces to the LA-ED model.

\subsubsection{Cost and reliability performance}
Table \ref{L6differentGammasMetrics} shows the performance of the two models: Column 2 for LA-ED, and Columns  3 to 7 for Rob-ED with different $\Gamma^w$'s. The best average total cost of the Rob-ED model is achieved at $\Gamma^w=0.5$, where the average cost of Rob-ED is $7.1\%$ lower than that of LA-ED; at the same time, Rob-ED is able to reduce the standard deviation of the cost by $41.2\%$. 
We can also see that as $\Gamma^w$ increases to $1.0$, the robust ED can reduce the std of cost by $82.1\%$, with the average cost reduced by $3.75\%$. The shortage event frequency of the robust ED model is decreased by up to $80.1\%$ and the associated penalty cost is reduced by $97.3\%$ at $\Gamma^w=1.0$. The change in penalty costs also implies that Rob-ED incurrs less amount of constraint violation than LA-ED, when penalty occurs. The results show that the robust ED model is effective at improving economic efficiency and reducing risk associated with the dispatch solution, where the risk exactly comes from the highly uncertain wind power. As will be shown in Section \ref{subsection:118bus}, more significant savings on cost and improvement over reliability are achieved for the 118-bus system.
\begin{table}[h!]
\renewcommand{\arraystretch}{1}
\caption{Performance of robust and deterministic ED}
\label{L6differentGammasMetrics}
\centering
\begin{tabular}{r|c|c|c|c|c|c}
\toprule
& LA-ED & \multicolumn{5}{c}{Rob-ED} \\
\hline
$\Gamma^w$ & {\bf 0.0} & 0.1 & 0.3 & {\bf 0.5} & 0.7 & {\bf 1.0} \\
\hline
Total Cost Avg (\$) & {\bf 771.1} & 758.5 & 734.0 & {\bf 716.0} & 718.2 & {\bf 742.2} \\
Total Cost Std (\$) & {\bf 1231} & 1172 & 1000 & {\bf 723} & 513 & {\bf 221} \\
Penalty Avg (\$) & {\bf 88.2} & 77.1 & 54.2 & {\bf 30.6} & 15.8 & {\bf 2.4} \\
Penalty Freq (\%) & {\bf 1.41} & 1.21 & 0.95 & {\bf 0.67} & 0.46 & {\bf 0.28} \\
\bottomrule
\end{tabular}
\end{table}

\subsubsection{Operational insights}\label{sec:OperInsight}
We also want to gain some insights about the operational characteristics of the robust model. Table \ref{L6differentGammasResults} shows average thermal generation (Therm avg) and wind generation (Wind avg) of the two models. 
We can see that the robust ED model on average tends to increase the use of thermal generation and curtail wind output: At $\Gamma^w=0.5$, Therm avg is up by $4.3\%$ and Wind avg down by $8.1\%$, comparing to LA-ED; at $\Gamma^w=1.0$, Therm avg is up by $16.1\%$ and Wind avg is down by $24.9\%$.

\begin{figure}[!t]
\centering
\includegraphics[width=3.5in]{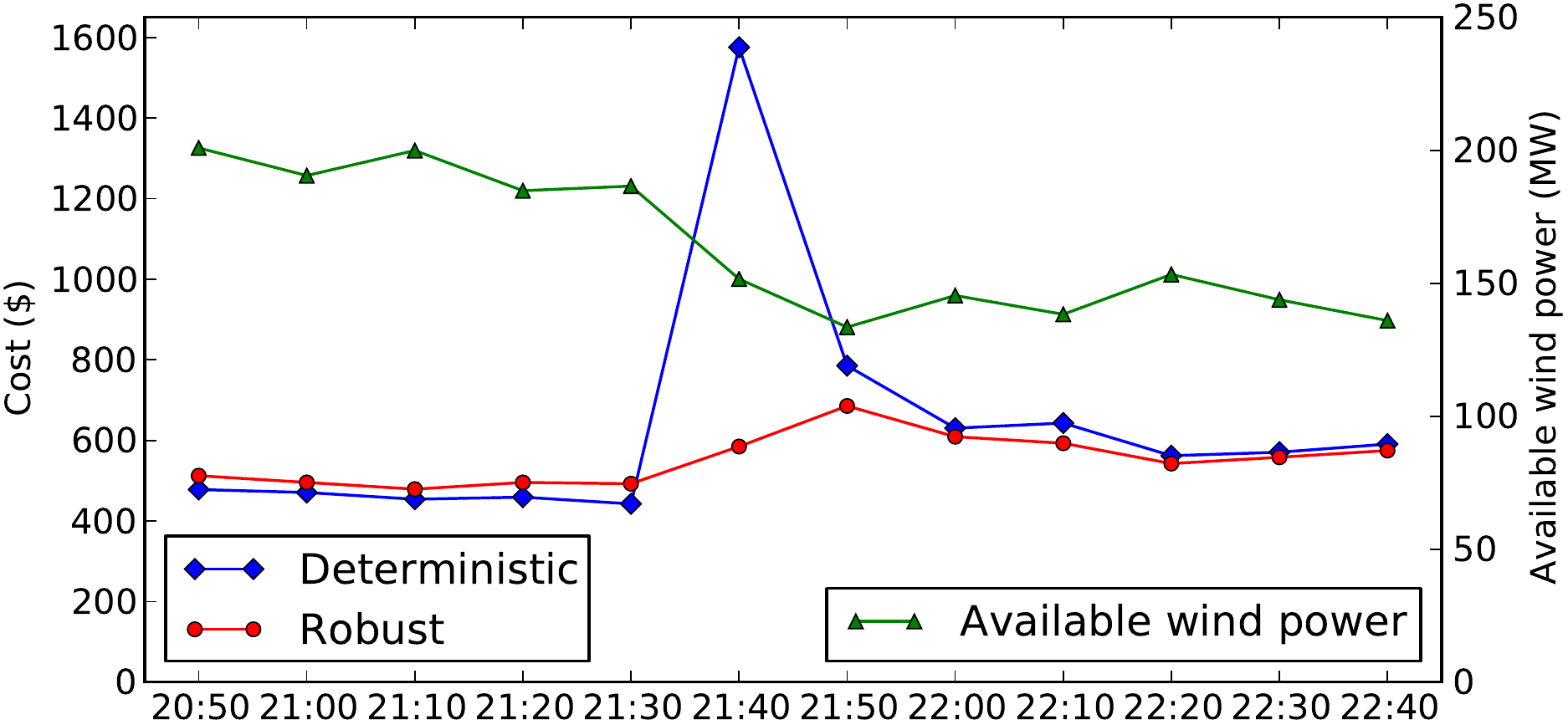}
\caption{A snapshot of the product cost of LA-ED and Rob-ED with $\Gamma^w=0.5$ when available wind power suddenly drops down.}
\label{CostL6Gamma00And05}
\end{figure}

Fig. \ref{CostL6Gamma00And05} shows a typical snapshot from simulation. Available wind power starts a fast and large drop at 21:30 (green curve), the deterministic LA-ED runs short of ramping capacity and incurs a spike of penalty cost (blue curve), while the system under robust ED is much less affected by this sudden wind event (red curve). The example shows that when the system has significant wind penetration, properly balancing wind and thermal generation becomes very important for system reliability.

The insight is the following. The two-stage robust ED computes wind scenarios over the future periods that are the most detrimental to the system, and makes the optimal dispatch solution to prepare the system against these scenarios. The worst-case wind scenarios often correspond to scenarios with large wind variation between periods as shown in Fig. \ref{CostL6Gamma00And05}. The robust ED model hedges against the potential large swing of wind by increasing thermal generation and moderately curtailing some wind output. In this way, the system maintains enough ramping capability to deal with potential sudden loss of available wind power. The balance between thermal and wind generation is controlled by the value of $\Gamma^w$ of the uncertainty sets as shown in Table \ref{L6differentGammasResults}. In other words, the robust ED determines the optimal ramping schedule of thermal generators, rather than resorting to prefixed operation rules.



\begin{table}[!t]
\renewcommand{\arraystretch}{1}
\caption{Operational Aspects of Robust and Deterministic ED}
\label{L6differentGammasResults}
\centering
\begin{tabular}{r|c|c|c|c|c|c}
\toprule
& {LA-ED} & \multicolumn{5}{c}{Rob-ED} \\
\hline
$\Gamma^w$ & {\bf 0.0} & 0.1 & 0.3 & {\bf 0.5} & 0.7 & {\bf 1.0} \\
\hline
Therm avg (MW) & {\bf 164.6} & 165.2 & 167.5 & {\bf 171.7} & 178.6 & {\bf 191.1} \\
Wind avg (MW) & {\bf 87.9} & 87.2 & 85.0 & {\bf 80.8} & 74.0 & {\bf 61.5} \\
\bottomrule
\end{tabular}
\end{table}

\subsubsection{Comparing to look-ahead ED with reserve}\label{sec:ReserveRules}
Reserve is an engineering approach to handle net load uncertainty in a deterministic ED model. Typically, when UC is solved, reserve levels for the next day are co-optimized, and later in real time operation, reserves are used in cases of unexpected net load variations and other contingencies.
Consider the following look-ahead ED model with reserve requirement (Res-LA-ED). The LA-ED model is complemented with reserve variables $R_{it} \in [0, \overline{R}_{it}]$, equations \eqref{LimitsG1} and \eqref{LimitsG2toT} are replaced by
\[\underline{p}_{it}^g \leq p^g_{it} \leq \overline{p}_{it}^g - R_{it} \quad \forall \, i \in \mathcal{N}^g, \, t=1,\dots,T,\]
and the following reserve requirement constraints are added:
\[\sum_{i \in \mathcal{N}^g} R_{it} \geq R^{req}_t \quad \forall \, t=1,\dots,T.\]
We test the performance of this model for different reserve requirement levels $R^{req}_t$. We select $R^{req}_t$ as a fraction of the total forecasted net load at time $t$ (i.e. forecast of total demand minus total available wind power), and modify this proportion, denoted as ``ResFactor'' \cite{papavasiliou2013multiarea}. Table \ref{ResLA-ED} presents the performance of Res-LA-ED under different values of ResFactor, as well as that of Rob-ED with $\Gamma^w = 0.5$.

From these results we can see that this reserve rule can improve the performance of LA-ED in both cost effectiveness and reliability, when the reserve requirement is properly chosen (ResFactor at $2.5\%$). As ResFactor increases, the reliability (Cost Std) keeps improving with the tradeoff of an increasing Avg Cost; the penalty cost and frequency are also reduced.

If we compare Res-LA-ED with Rob-ED, we can observe that the performance of Rob-ED is significantly better than the best of the three Res-LA-ED  test cases: the Cost Avg is reduced by at least $7.14\%$ (against $\text{ResFactor}=2.5\%$); the Cost Std is improved by at least $37.4\%$ (against $\text{ResFactor}=10\%$); the penalty cost is reduced by at least $57.2\%$, and the penalty frequency is reduced by at least $50.3\%$ (both against $\text{ResFactor}=10\%$).

\begin{table}[!t]
\renewcommand{\arraystretch}{1}
\caption{Performance of look-ahead ED with reserve}
\label{ResLA-ED}
\centering
\begin{tabular}{c|c|c|c|c|c}
\toprule
& LA-ED & \multicolumn{3}{|c|}{Res-LA-ED} & Rob-ED \\
\hline
ResFactor (\%) & 0 & 2.5 & 5 & 10 & $\Gamma^w=0.5$ \\
\hline
Cost Avg (\$) & 771.1 & 770.0 & 773.3 & 790.3 & 716.0 \\
Cost Std (\$) & 1231 & 1223.8 & 1211.8 & 1155.1 & 723 \\
Penalty Avg (\$) & 88.2 & 86.7 & 84.8 & 71.6 & 30.6 \\
Penalty Freq (\%) & 1.41 & 1.45 & 1.69 & 1.35 & 0.67 \\
\bottomrule
\end{tabular}
\end{table}

\subsection{Dynamic uncertainty sets versus static uncertainty sets} \label{subsection:staticvsdynamic}
In this section, we compare adaptive robust ED equipped with dynamic uncertainty sets with the same robust ED model using static uncertainty sets. The goal is to study the benefits of dynamic uncertainty sets for modeling dynamic relations of wind power uncertainty across time stages and spatial locations.

We use dynamic uncertainty sets \eqref{eq:windspeed_dus} with $L=6$ as before (denoted as ``DUS''), and construct two static uncertainty sets: one ignores the temporal correlation in \eqref{eq:windspeed_dus} (denoted as ``SUS1''), the other further ignores spatial correlations (denoted as ``SUS2''). Note that both SUS1 and SUS2 are special cases of the dynamic uncertainty sets for $L=0$, i.e. the uncertainty sets at different time intervals are independent of each other. To have a fair comparison, both in SUS1 and SUS2, $\mb{g}_t$ is improved after estimating $\mb{B}$ to force a persistent forecast of wind speeds for the nominal trajectory (improving the accuracy of the nominal trajectory considered).

\begin{figure}[!t]
\centering
\includegraphics[scale=0.45]{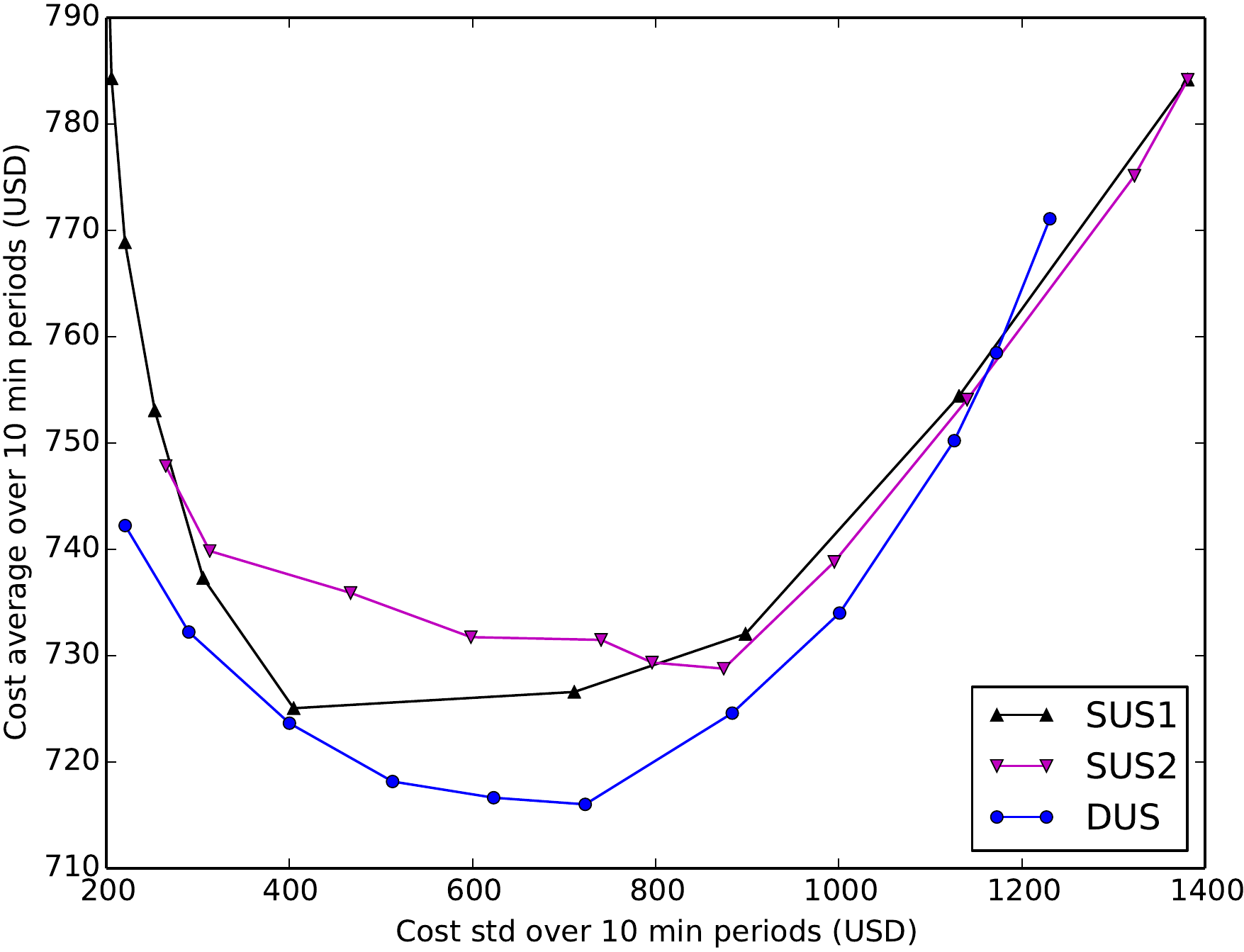}
\caption{Cost std and cost average obtained for the policies determined by the different models with $\Gamma^w = 0.0, 0.1,..., 1.0$}
\label{TenMinStdDevVsCost}
\end{figure}

Fig. \ref{TenMinStdDevVsCost} plots the standard deviation of the cost per 10 minutes interval (x-axis) versus the average of this cost (y-axis) for DUS, SUS1 and SUS2 with different values of $\Gamma^w$.
On each curve, the right most point corresponds to $\Gamma^w=0$, i.e. the deterministic LA-ED model. As $\Gamma^w$ increases, both the average and std of the cost start to decrease, then after a certain apex value of $\Gamma^w$ around $0.4$ to $0.5$, the std keeps decreasing but the average cost starts to increase. This behavior endows a ``U'' shape for all three curves. Every point on the right half of the ``U" shape for $\Gamma^w$ smaller than the apex value can be strictly improved in both average and std of cost by increasing $\Gamma^w$, while every point on the left half of the ``U'' shape cannot be strictly improved without trading off between average and std of the cost. In other words, on the right half of the curve, each point is {dominated} by the points to its left, whereas on the left half, no point is dominated by any other. Therefore, the left part of each curve shows the \emph{Pareto frontier} of cost average vs cost standard deviation performance of the associated robust ED model. The system should be operated on the Pareto frontier. This provides an informative guideline for choosing a proper $\Gamma^w$.


Comparing the Pareto frontiers of the three uncertainty sets in Fig. \ref{TenMinStdDevVsCost}, we can see that the dynamic uncertainty set has the lowest Pareto frontier, which means that to retain a same level of average cost, the robust ED with dynamic uncertainty sets achieves the lowest std (i.e. the highest reliability); or, to maintain a same level of std (i.e. reliability), the robust ED with dynamic uncertainty sets incurs the lowest cost. That is, robust ED with DUS \emph{dominates} robust ED with static uncertainty sets. Between the two static uncertainty sets, SUS1 (that considers spatial correlation) dominates SUS2, which has neither temporal nor spatial correlation.



The static uncertainty set SUS2 is the first budgeted uncertainty set proposed in the literature \cite{Bertsimas2004price} and has inspired its application in modeling net load uncertainty \cite{Sun2013robustUC}. Works in \cite{Guan2012pumphydro,Zeng2012robustUC} further introduced budget constraints over time periods to limit the total variations of uncertain demand over the entire or part of the planning horizon. Now, we compare these static uncertainty sets with additional time budgets with DUS. It is worth emphasizing that the  fundamental difference between DUS and SUS remains the same for DUS and SUS with time budgets.

We modify the uncertainty sets SUS1 and SUS2 with the following time budget constraint:
\[\sum_{t=2}^T \sum_{i \in \mathcal{N}^w} |u_{it}| \leq \Gamma^T \Gamma^w \sqrt{N^w} \sqrt{T-1},\]
where $T=9$ is the number of periods in the multi-period Rob-ED, and time budget parameter $\Gamma^T = 0.5, 1, 2$. Note that static uncertainty sets without time budget are equivalent to one with a large time budget as when $\Gamma^T \geq \sqrt{8}$ the time budget constraint becomes redundant.

\begin{figure}[!t]
\centering
\includegraphics[scale=0.45]{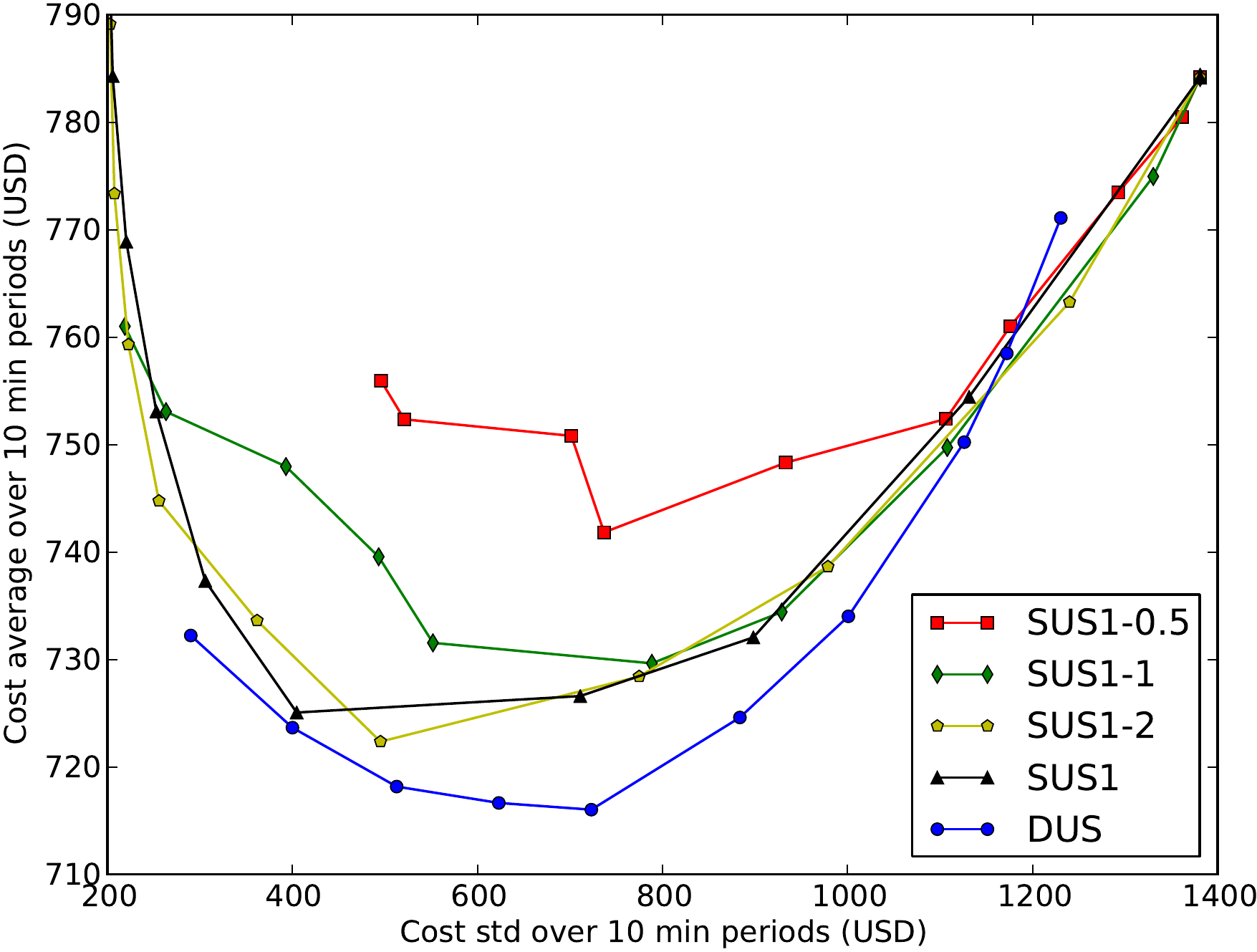}
\caption{Cost std and cost average obtained for the policies determined by DUS and SUS1, with $\Gamma^T = 0.5, 1, 2$ for SUS1 and with $\Gamma^w = 0.0, 0.1,..., 1.0$.}
\label{TenMinStdDevVsCostL0BfullTimeBudget}
\end{figure}

Fig. \ref{TenMinStdDevVsCostL0BfullTimeBudget} plots the std of cost per 10 min interval (x-axis) versus the average of that cost (y-axis) for DUS and SUS1 with additional time budgets. The curve denoted by SUS1-0.5 means the SUS1 uncertainty set with time budget $\Gamma^T=0.5$ and $\Gamma^w$ varies from $0.0$ to $1.0$.  Among the three curves based on SUS1 with time budgets, we can see that Rob-ED achieves a better Pareto frontier for higher values of time budget (the red curve for SUS1-0.5 is dominated by the green curve for SUS1-1, which is further dominated by SUS1-2). SUS1 without time budget (or equivalently with a time budget $\Gamma^T\geq\sqrt{8}$) has a frontier comparable to the SUS1-2. Furthermore, all four SUS1 based curves are clearly dominated by the DUS curve.

\begin{figure}[!t]
\centering
\includegraphics[scale=0.45]{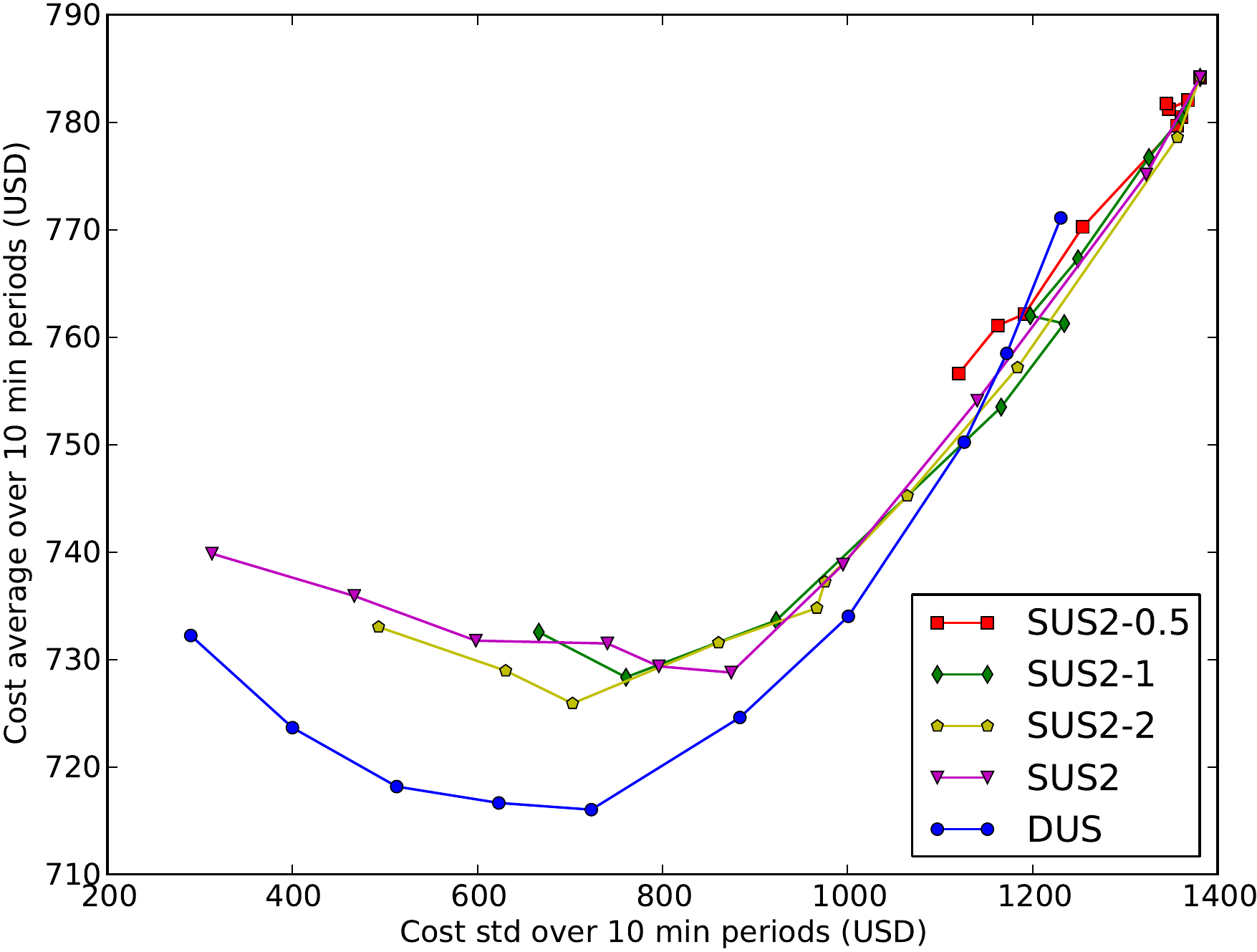}
\caption{Cost std and cost average obtained for the policies determined by DUS and SUS2, with $\Gamma^T = 0.5, 1, 2$ for SUS2 and with $\Gamma^w = 0.0, 0.1,..., 1.0$ for all policies}
\label{TenMinStdDevVsCostL0BdiagTimeBudget}
\end{figure}

Fig. \ref{TenMinStdDevVsCostL0BdiagTimeBudget} presents a similar comparison for SUS2 with time budgets. Here, the dominance of DUS over static uncertainty sets with time budgets is more eminent.

\subsection{Impact of system ramping capacity} \label{subsection:sensitivities}
\begin{figure}[!t]
\centering
\includegraphics[scale=0.45]{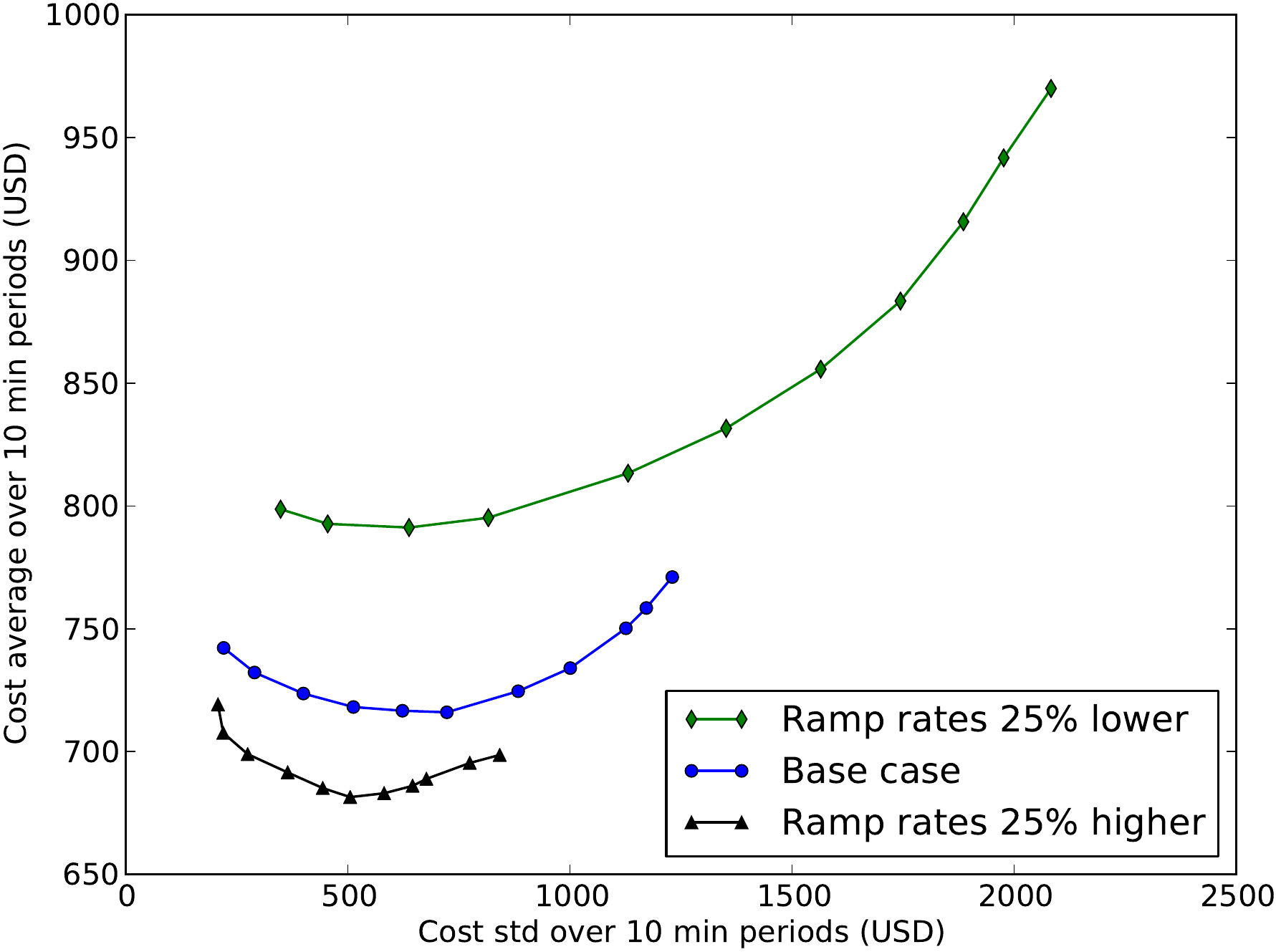}
\caption{Cost std and cost average obtained for the Rob-ED with DUS for $\Gamma^w = 0.0, 0.1,..., 1.0$, under modified ramping rates.}
\label{TenMinStdDevVsCostForL6Changes}
\end{figure}

In this section, we study the relationship between system ramping capacity and the performance of robust ED models. The intuition is that higher ramping rates better prepare the system to deal with high variation of wind output. We want to see how much benefit the robust ED model provides under different system ramping capacities. Fig. \ref{TenMinStdDevVsCostForL6Changes} summarizes the computational results for three scenarios: base case with no change in ramping rates,  and $-25\%$ or $+25\%$ change on each generator's ramping rates.

We can see that the robust ED model saves the average cost by $7.1\%$ in the base case (the same numbers as in Section \ref{subsection:robustvsdeterministic}) comparing with the look-ahead ED; the saving increases to $21.2\%$ for the reduced ramping case; even for the system with $25\%$ more ramping for every generator, the robust ED still demonstrates a $3.7\%$ saving in average cost over LA-ED. This demonstrates the clear benefit of Rob-ED over a wide range of system ramping conditions.

\subsection{Considering both demand and wind uncertainty} \label{subsection:demanduncertainty}
In this section, we further incorporate traditional demand uncertainty into the robust ED model, using the static uncertainty sets \eqref{eq:sus}, where $\overline{d}_{jt}$ and $\hat{d}_{jt}$ are selected as the mean and std of demand from previously realized values. The parameter $\Gamma^d$ limits the total deviation of demand from its forecast. In simulation, the demand $d_{jt}$ of each load $j$ at each time period $t$ is independently generated as a normal random variable with a std that equals a 5\% of its mean, and is truncated to be nonnegative. Therefore, the generated random demand can be outside the uncertainty set. The choice of $\Gamma^d$ controls the size of the demand uncertainty set.

\begin{figure}[!t]
\centering
\includegraphics[scale=0.45]{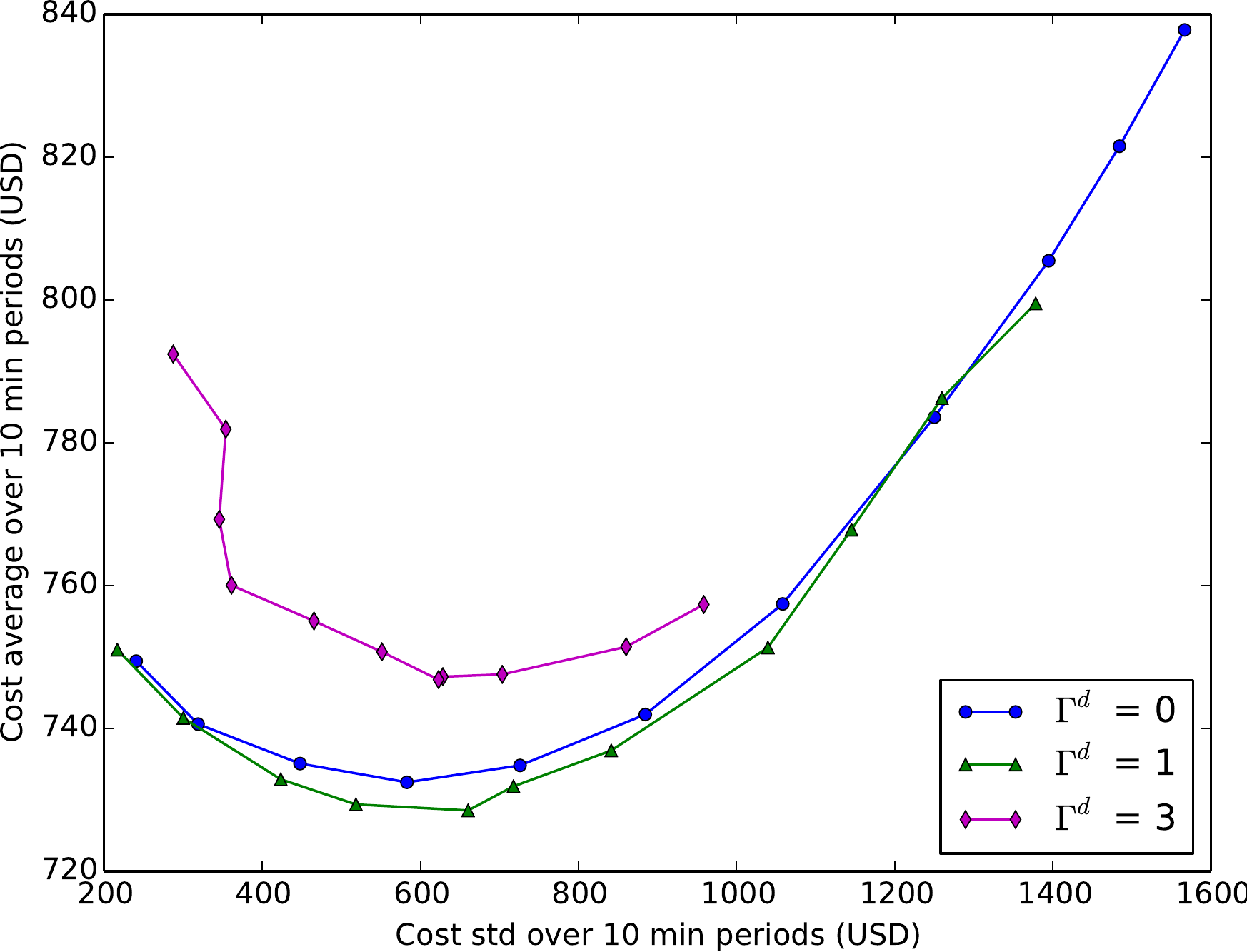}
\caption{Cost std and average obtained by the Rob-ED with $\Gamma^w = 0.0,...,1.0$ and $\Gamma^d=0,1,3$.}
\label{Demandstd5percent-L6-TenMinStdDevVsCost}
\end{figure}

Fig. \ref{Demandstd5percent-L6-TenMinStdDevVsCost} presents the performance of Rob-ED with dynamic uncertainty set for wind and static uncertainty sets \eqref{eq:sus} for load, at different values of $\Gamma^d,\Gamma^w$. At $\Gamma^d=0$, the uncertainty set for demand is a singleton containing the forecast value, i.e. only wind uncertainty is considered (blue curve). By considering an uncertainty set for load with $\Gamma^d=1$, the cost-reliability curve is shifted downward to the green curve, which consistently dominates the blue curve. The two curves are quite close though, which shows that wind is the dominating factor of uncertainty; the dynamic uncertainty set for wind significantly improves the system performance, while further incorporating load uncertainty improves the performance modestly. The purple curve for $\Gamma^d=3$ shows that too much conservatism in the load uncertainty model leads to inferior solutions. It again demonstrates that properly choosing the level of conservativeness of the uncertainty sets is critical to getting the best performance of the robust ED model. In particular, the best robust ED policy obtained by setting $\Gamma^d = 1,\Gamma^w = 0.6$ reduces the average cost by $13.1\%$ lower than that of the deterministic LA-ED with $\Gamma^d = \Gamma^w = 0$, and reduces the std of the cost by $58.1\%$. This makes the robust ED model very attractive.


\subsection{Performance of the alternating direction method for solving the second-stage problem} \label{subsection:alternatingmethodperformance}
As discussed in Section \ref{Subsection:SolutionMethod}, the proposed algorithm requires solving a bilinear program \eqref{eq:Qafterdual} in each iteration of the outer master problem. Therefore, to practically tackle large-scale problems, a fast and reliable method for the inner problem is needed. An alternating direction (AD) method is proposed in Section \ref{Subsection:SolutionMethod} for this purpose. In the literature, several exact MIP methods are proposed to solve the second-stage bilinear program (e.g. \cite{Guan2012pumphydro,Zeng2012robustUC,Zeng2013ccg}.) The MIP methods in \cite{Guan2012pumphydro,Zeng2012robustUC} rely on the special structure of the uncertainty sets used in their models, which are not shared by the dynamic uncertainty sets. The exact MIP method proposed in \cite{Zeng2013ccg} is based on the KKT conditions, which are applicable to general polyhedral uncertainty sets. Thus, we compare the AD algorithm to this MIP method.

In the experiment, we run the Rob-ED model in the rolling-horizon simulator for a 5-day horizon. This involves solving 720 Rob-ED models  of the form $\eqref{Formulation}$, which amounts to 1529 inner bilinear programs \eqref{eq:Qafterdual}. Every time, the bilinear program is solved by both the AD algorithm and the exact MIP method. We compare both running times and solution qualities.

The AD algorithm achieves convergence for all 1529 instances, and the average running time is 0.12s. The MIP method achieves convergence in 257 instances with an average time of 13.28s; for the remaining 1272 instances ($83.2\%$ of the total instances), the MIP method does not converge after 60s, and at that point the solution quality is still worse than the AD solutions (the objective value is on average $1.02\%$ worse than the AD solutions). Those MIP instances exceeding 60s do not achieve much improvement after running for another 10 min. In terms of solution quality, the AD solutions on average obtain an optimality gap of $3.73\%$ compared to the global optimum of the MIP solutions when MIP converges. These comparisons show that the AD algorithm is an effective and efficient heuristic for solving the bilinear program.



\subsection{Tests on 118-bus system} \label{subsection:118bus}
Extensive simulation is also conducted on the 118-bus system. The results for this larger system support similar conclusions as shown in the 14-bus system. The 118-bus system has $54$ generators of total $7220$ MW generation capacity and 273.2 MW/10min system ramping capacity. There are 186 lines with flow limits ranging between 280 MW and 1000 MW, and 91 loads. Total electricity demand is between 2485.7 MW (3:30 am) and 5982.9 MW (8:20 pm) with an average of 4735.0 MW. There are 8 wind farms, each with a capacity $p^{w,max} = 750$ MW. The average total available wind power at any time is 1882.7 MW, equivalent to $31.5\%$ of the peak demand. All the wind speeds used in simulation are real data collected from \cite{ceerewebsite}. Each robust ED model can still be quickly solved in about 20 seconds in the laptop described before. For the simulation of 35-day rolling horizon with a 10-min interval, we use a computer cluster \cite{isyecluster}.

Table \ref{L6differentGammasCase118} shows the performance of the deterministic LA-ED and the Rob-ED with dynamic uncertainty sets of lags $L=6$. From the table, we have the following observations:
\begin{enumerate}
\item[(1)] Rob-ED reduces the average cost by $43.4\%$ ($(15061-8528)/15061$) at a properly chosen $\Gamma^w=1.5$.
\item[(2)] Cost std is reduced by $87.7\%$ at $\Gamma^w=1.5$ and by $93.9\%$ at $\Gamma^w=2.0$.
\item[(3)] The average penalty cost is reduced by $98.4\%$ or $60.7$ times at $\Gamma^w=1.5$ and is almost eliminated at $\Gamma^w=2.0$. The frequency of penalty is $7.70\%$ by LA-ED, and is reduced to $0.12\%$ and $0.02\%$ by Rob-ED at $\Gamma^w=1.5$ and $2.0$, respectively.
\item[(4)] Rob-ED dispatches more thermal and curtails more wind. On average, the thermal generation is up by $12.7\%$ and $18.9\%$, and the wind generation is down by $24\%$ and $38.9\%$, at $\Gamma^w=1.5$ and $2.0$, respectively. This can be explained by a similar reasoning given in Section \ref{sec:OperInsight}, namely that the robust ED dispatches the thermal generation anticipating to a potential large drop of wind in the future, optimally balancing thermal and wind generation in the system.
\end{enumerate}


Comparing to the 14-bus system, the above results for the 118-bus system show a more significant benefit of the proposed Rob-ED model: the average operating cost is cut to almost half of the look-ahead ED, the cost variability is reduced by an order of magnitude, and the shortage events and penalty cost are almost eliminated. Table \ref{L6differentGammasCase118} also shows a Pareto frontier exists for the range of $\Gamma^w$ between $1.5$ to $2.0$.

\begin{table}[!t]
\renewcommand{\arraystretch}{1}
\caption{Performance of LA-ED and Rob-ED for 118-bus system}
\label{L6differentGammasCase118}
\centering
\begin{tabular}{r|c|c|c|c|c}
\toprule
& LA-ED & \multicolumn{4}{c}{Rob-ED} \\
\hline
$\Gamma^w$ & 0.0 & 0.5 & 1.0 & 1.5 & 2.0 \\
\hline
Cost Avg (\$) & 15061 & 12193 & 8914 & 8528 & 9075 \\
Cost Std (\$) & 38138 & 30903 & 14671 & 4703 & 2325 \\
Penalty Avg (\$) & 7775 & 4835 & 1214 & 126 & 1 \\
Penalty Freq (\%) & 7.70 & 4.74 & 1.45 & 0.12 & 0.02 \\
Therm Avg (MW) & 2969 & 3007 & 3132 & 3399 & 3660 \\
Wind Avg (MW) & 1758 & 1723 & 1602 & 1336 & 1075 \\
\bottomrule
\end{tabular}
\end{table}

\section{Conclusion} \label{Section:Conclusion}
In this paper, we present an adaptive multi-period robust ED model and dynamic uncertainty sets for power system economic dispatch under high penetration levels of wind resources. The adaptive multi-period robust ED model mimics the physical dispatch procedure by using a two-stage decision making structure and a rolling-horizon framework. Dynamic uncertainty sets explicitly model the relationship between uncertainties across decision stages and capture the temporal and spatial correlations of wind power output in multiple wind farms: the proposed dynamic uncertainty sets with linear dynamics in this paper have general and computationally tractable structure; and the proposed data-driven estimation procedures are easy to implement. We also develop a simulation platform that integrates the optimization engine and data analysis tools for updating uncertainty sets.

Extensive simulation using real wind data shows that the proposed robust ED framework outperforms look-ahead ED models with and without reserves which recently attracted considerable interests in practice, and robust ED models with static uncertainty sets. Both cost efficiency and system reliability are substantially improved. Also, the robust ED model gives an entire Pareto frontier of operating cost and reliability, which provides an informative guideline for choosing uncertainty set parameters and system operating points. The proposed robust ED model and dynamic uncertainty sets are flexible enough to incorporate several extensions, such as using transformed wind speeds, bids with piecewise linear costs, and including other types of uncertain renewable energy sources.




\appendices



\ifCLASSOPTIONcaptionsoff
  \newpage
\fi



\bibliographystyle{IEEEtranS}
\bibliography{WindDispatchBib}

\end{document}